\documentclass[J]{kkmjart}

\usepackage{mathrsfs}
\usepackage{amsthm, amscd}
\usepackage{amsfonts,amsmath,amssymb,amsthm,cite,enumerate,epsfig,fancyhdr,subfigure}

\usepackage[hidelinks]{hyperref}
\hypersetup{
    colorlinks=false,
    pdfborder={0 0 0},
}

\newtheorem{theorem}{Theorem} 
\newtheorem{lem}{Lemma}

\newtheorem{app}{Application}

\newtheorem{corollary}{Corollary}

\newtheorem{f}{Fact}

\theoremstyle{definition}

\newtheorem{example}{Example}

\begin{document}

\title[THE PROPERTIES OF BANACH ALGEBRAS AFFECTED BY DIFFERENTIAL IDENTITIES
 ]{ THE PROPERTIES OF BANACH ALGEBRAS AFFECTED BY DIFFERENTIAL IDENTITIES}

\author[M. MOUMEN]{Mohamed MOUMEN }
\address{{\bf National School of Applied Sciences}
\\Department of Mathematics
\\Ibn Zohr University
\\ Agadir, Morocco }
\email{mohamed.moumen@edu.uiz.ac.ma}

\author[L. TAOUFIQ]{Lahcen TAOUFIQ} 
\address{{\bf National School of Applied Sciences}
\\Department of Mathematics
\\Ibn Zohr University
\\ Agadir, Morocco}
\email{l.taoufiq@uiz.ac.ma}

%\thanks{$\dagger$ This work was supported by KKMS.}
%\thanks{$*$ Corresponding author.}

%%%%%%%%%%%%%%%%%%%%%%%%%%%%%%%%%%%%%%%
\begin{abstract}
Let $\mathcal{B}$ be a Banach algebra. The interest of this article lies in the study of the commutativity of $\mathcal{B}$ if certain specific algebraic identities hold over a non-empty open subset of $\mathcal{B}$. The limitations imposed in the
hypothesis of our results are justified by some examples.
\end{abstract}
%%%%%%%%%%%%%%%%%%%%%%%%%%%%%%%%%%%%%%%
\maketitle

%%%%%%%%%%%%%%%%%%%%%%%%%%%%%%%%%%%%%%%
\section{Introduction}
%==================================================

%\numberwithin{equation}
%\renewcommand{\thefootnote}{\*}

%========================================

	\hspace{0.5cm} 
	$\mathcal{B}$ will  denote a Banach algebra with center $Z(\mathcal{B})$. If for every $x, y \in \mathcal{B},$ $x\mathcal{B}y = 0$ means that either $x = 0~~\text{or}~~y=0$, then $\mathcal B $ is said to be prime. 
	The anti-commutator $xy + yx$ will be represented by the symbol $x \circ y$, whereas the commutator $[x, y]$ will stand for the Lie product $xy-yx$. An additive mapping $d$ in $\mathcal{B}$  is said to be a derivation if $d(xy) = d(x)y + xd(y)$ for all $x, y \in \mathcal{B}.$ We recall that $f: \mathcal{B} \longrightarrow  \mathcal{B}$ is called a homomorphism if $f(xy)=f(x)f(y)$ and $f(x+y)=f(x)+f(y)$ holds for any $x, y \in \mathcal{B}.$ %Following Bre$\breve{s}$ar \cite{bresar}, 
	An additive mapping $F: \mathcal{B} \longrightarrow  \mathcal{B}$ is said to be a generalized derivation associated with a derivation $d$ if $F (xy) = F (x)y + xd(y)$ for all $x,y\in\mathcal{B}$. Basic examples are derivations and generalized inner derivations, i.e., maps of type $x \mapsto ax+xb $ for some $a, b\in \mathcal{B}$. A mapping $f:\mathcal{B}\longrightarrow \mathcal{B}$ is  centralizing on a part $S$ of $\mathcal{B}$ if $[f(x),x]\in Z(\mathcal{B})$ for all $x\in S$, in particular if $[f(x),x]=0$ for all $x\in S$ then $f$ is said to be commuting on $S$ .\\ In 1955, Divinsky \cite{divinsky} established that a simple Artinian ring is commutative if it possesses a commuting non-trivial automorphism, which began the development of commuting and centralizing mappings.
	Posner\cite{posner} demonstrated two years later that a prime ring with a nonzero centralizing derivation must be commutative. 
	
	Over the last few decades, several authors have proved commutativity theorems for prime and semiprime rings admitting automorphisms or derivations that are commuting
    or centralizing on an appropriate subset of the ring
	(see \cite{r7} for a partial bibliography). In the direction of Banach algebras, Yood  \cite{Yood94}  showed that a semiprime complex Banach algebra $\mathcal{B}$ must be commutative if there are non-empty open subsets $\mathcal H_1$ and  $\mathcal H_2$
of $\mathcal{B}$ such that for each $(x,y) \in \mathcal H_1\times \mathcal H_2$  there are positive integers
$n = n(x, y)>1$ and $\;m = m(x, y)>1$ depending on $x$ and $y$ such
that either $[x^n,y^m]\in Z(\mathcal{B})$ or $x^ny^m \in Z(\mathcal{B}).$ Recently,  M. Moumen, L. Taoufiq and L. Oukhetite \cite{Moumen}  proved that  a  prime Banach algebra $\mathcal{B}$ must be commutative if there are non-empty open subsets $\mathcal{H}_{1},\mathcal{H}_{2}$  and  a non-injective continuous derivation $d$  satisfying for all $(x,y)\in \mathcal{H}_{1}\times \mathcal{H}_{2}$ there are strictly  positive integers $n=n(x,y),m=m(x,y)$ such that $d(x^{n}y^{m})+ [x^{n},y^{m}]\in Z(\mathcal{B}).$ In \cite{Mou} M. Moumen, L. Taoufiq and A. Boua proved similar results
(see \cite{Mou} and \cite{mm} for additional details).\\

The importance of commutativity in many disciplines such as spectral theory, information security, and especially decryption theory, is the reason why we chose to continue in this domain.
    %So the earlier works provided as the inspiration for this article, we continue to investigate the commutativity of Banach algebra by employing automorphism.
    Our method is based on some functional analysis theorems as  Baire's category theorem. Among the results, we have demonstrated that a prime Banach algebra $\mathcal{B}$ must be commutative if it has non-empty open subsets$\mathcal{H}_{1},\mathcal{H}_{2}$ and a continuous automorphism $f$ satisfying for any $(x,y)\in \mathcal{H}_{1}\times \mathcal{H}_{2}$ there exist two strictly positive integers $p=p(x,y), q=q(x,y)$ such that $f(x^p y^q) + x^p \circ y^q \in Z(\mathcal{B})$. Further, we provide applications for our results. \\
    %The article is concluded with a number of examples.\\
    
    Throughout this article, we shall make some use of the following well-known results (see \cite{Mou} for further details).
\begin{f}\label{rm}
 Let $\mathcal{R}$ be a prime ring.\\
1. If $x\in Z(\mathcal{R})$ and $xy\in Z(\mathcal{R}),$ then $x=0$ or $y\in Z(\mathcal{R}).$\\
%2.  $Z(\mathcal{R})$ does not admit any zero divisors.\\
2. If  $\mathcal{R}$ has a nonzero centralizing derivation, then $\mathcal{R}$ is commutative.
\end{f}
\section{Automorphism in Banach algebras}
The development of the confirmation of our main ideas depends on the following lemma, attributed to  Duncan and  Bonsall {\rm \cite{Bonsal}}.
\begin{lem}\label{4}
Let $s(t)=\displaystyle\sum_{i=0}^{n} t^{i}b_{i} $ be a polynomial with coefficients in $\mathcal{B}$ (complex or real Banach algebra). If for an infinite set of real values of $t,\;s (t) \in C,$
where $C$ is a closed linear subspace of $\mathcal{B},$ then every $b_i$ lies in $C.$

\end{lem}
\begin{theorem}\label{a}
Let $ \mathcal{B} $ be a prime Banach algebra, $\mathcal{H}_{1}$ and $\mathcal{H}_{2}$ are non-void open subsets of $ \mathcal{B} $.  $ \mathcal{B} $ must be commutative if admits a continuous automorphism $f$ satisfying: $f(x^{p}y^{q})+x^{p}\circ y^{q}\in Z(\mathcal{B})$ for all $(x,y)\in \mathcal{H}_{1}\times \mathcal{H}_{2}$ where $p, q$ are not  fixed but  depends on the pair of elements $x$ and $y$.
\end{theorem}
\begin{proof}
We define the following set for any $(p,q)\in\mathbb N^{*}\times \mathbb N^{*}$:
 $$O_{p,q}=\{(x,y)\in \mathcal{B}^{2}\,|\, f(x^{p}y^{q})+ x^{p}\circ y^{q}\notin Z(\mathcal{B} )\}.$$
We note that $(\cap O_{p,q})\cap (\mathcal H_{1} \times \mathcal H_{2})= \emptyset$, indeed:\\
If they is $(a,b)\in\mathcal H_{1}\times \mathcal H_{2} $ such that $(a,b)\in O_{p,q}$ for all $(p,q)\in\mathbb N^{*}\times \mathbb N^{*},$ then \\$f(x^{p}y^{q})+ x^{p}\circ y^{q}\notin Z(\mathcal{B})$ for all $(p,q)\in\mathbb N^{*}\times \mathbb N^{*},$ this contradicts the theorem's assumption.\\
We now verify that every $O_{p,q}$ is open in $\mathcal{B}\times\mathcal{B}$. In other words, we demonstrate that the complement of $O_{p,q}$, $O_{p,q}^c$, is closed. In fact, to do this, we consider  a sequence
$((x_{k},y_{k}))_{k\in\mathbb{N}}
\subset O^c_{p,q}$ that converges to $(x,y)\in\mathcal{B}\times\mathcal{B}$.
Since $((x_{k},y_{k}))_{k\in\mathbb{N}}\subset O^c_{p,q}$, so
\begin{center}
$f((x_{k})^{p}(y_{k})^{q})+ (x_{k})^{p}\circ (y_{k})^{q}\in Z(\mathcal{B})$ for all $k\in\mathbb{N}.$
\end{center}
Given that $f$ is a continuous maps, the sequence  $(f((x_{k})^{p}(y_{k})^{q})+ (x_{k})^{p}\circ (y_{k})^{q})_{k\in\mathbb{N}}$ converges to
$ f(x^{p}y^{q})+ x^{p}\circ y^{q} $,
know that  $Z(\mathcal{B})$ is closed, we obtain $ f(x^{p}y^{q})+ x^{p}\circ y^{q}\in Z(\mathcal{B}).$
As a result,  $ (x,y)\in O^c_{p,q}$ and  $O^c_{p,q}$ is closed(i.e $O_{p,q}$ is open).\\
Whenever $ O_{p,q} $ is dense, the Baire's category theorem informs us that their intersection is also dense, which conflicts with  $(\cap O_{p,q})\cap (\mathcal H_{1} \times \mathcal H_{2})= \emptyset$. Hence, there is  $(n,m)\in\mathbb N^{*}\times \mathbb N^{*}$ such that  $ O_{n,m} $ is not a dense open and there exists a non-void open subset $\mathcal {O} \times \mathcal O^{'}$ in $O^c_{p,q}$ such that:
\begin{center}
$f(x^{n}y^{m})+ x^{n}\circ y^{m}\in Z (\mathcal{B})\,\,\,\mbox{for all}\,\,\,x\in \mathcal O, y\in \mathcal O^{'},$
\end{center}
Fix $ y\in \mathcal O^{'}.$ Let $x \in \mathcal O $ and $z\in\mathcal{B}$, then $x+tz\in \mathcal O$ for all sufficiently small real $t$. Therefore  $P(t)= f((x+tz)^{n}y^{m})+(x+tz)^{n}\circ y^{m}\in Z(\mathcal{B})$. We are capable of writing
$$P(t)=B_{n,0}(x,z,y)+tB_{n-1,1}(x,z,y) + t^2B_{n-2,2}(x,z,y)+...+t^nB_{0,n}(x,z,y).$$
By Lemma $\ref{4},$ we conclude that $B_{0,n}(x,z,y)=f(z^{n}y^{m})+z^{n}\circ y^{m}\in Z(\mathcal{B}).$ Accordingly, chosen to give $ y\in \mathcal O^{'}$ we have
$ f(x^{n}y^{m})+x^{n}\circ y^{m}\in Z(\mathcal{B})$ for all $x\in \mathcal{B}.$ We reverse the roles of $\mathcal O $  and $\mathcal O^{'} $ in the above settings with  $x$ is fixed in $\mathcal{B}$, we find that  $f(x^{n}y^{m})+x^{n}\circ y^{m}\in Z(\mathcal{B})$ for all $(x,y)\in \mathcal B^{2}.$ \\
Replacing $x$ by $x^{m}$ and $y$ by $y^{n}$, we find that
$$
f(x^{nm}y^{mn})+ x^{nm}\circ y^{mn}\in Z(\mathcal{B})\,\,\,\mbox{for all}\,\,\, x,y\in\mathcal{B}.
$$
By changing $x$ and $y$, we obtain
$$
f(y^{nm}x^{mn})+ y^{nm}\circ x^{mn}\in Z(\mathcal{B})\,\,\,\mbox{for all}\,\,\, x,y\in\mathcal{B}.
$$
Then $$(f(x^{nm}y^{mn})+ x^{nm}\circ y^{mn})-(f(y^{nm}x^{mn})+ y^{nm}\circ x^{mn})\in Z(\mathcal{B})\,\,\,\mbox{for all}\,\,\, x,y\in\mathcal{B}.$$
Since  $$y^{nm}\circ x^{mn}=x^{nm}\circ y^{mn}\,\,\,\mbox{for all}\,\,\, x,y\in\mathcal{B},$$
we obtain
$$
f([x^{nm},y^{mn}])\in Z(\mathcal{B})\,\,\,\mbox{for all}\,\,\, x,y\in\mathcal{B}.
$$
According to Theorem $2.2$ of \cite{mmm}, we come to the conclusion that $\mathcal{B}$ is commutative.
\end{proof}
\begin{theorem}\label{b}
Let $\mathcal{H}_{1}$ and $\mathcal{H}_{2}$ are non-void open subsets of a prime Banach algebra $\mathcal{B}.$
 If  a continuous automorphism $f$ of $\mathcal{B}$, satisfying  the following condition:
 $$ (\forall (x,y)\in \mathcal{H}_{1}\times \mathcal{H}_{2})(\exists (p,q)\in\mathbb N^{*}\times \mathbb N^{*}\;\; \mbox{such\; \;  that}\; \;
  f(x^{p}\circ y^{q})+[x^{p},y^{q}]\in Z(\mathcal{B}),$$
 then $ \mathcal{B}$ is  commutative.
\end{theorem}
\begin{proof}
In the same way as Theorem $\ref{a}$'s proof, we arrive at the conclusion that there is $(n,m)\in\mathbb N^{*}\times \mathbb N^{*}$ such that
$$
f(x^{n}\circ y^{m})+ [x^{n},y^{m}]\in Z(\mathcal{B})\:\:\mbox{for all}\:\: x, y\in\mathcal{B}.
$$
If $ y=x,$ then $ 2f(x^{n+m})\in Z(\mathcal{B})$ for all $ x\in\mathcal{B},$ since $\mathcal{B}$ is 2-torsion free  therefore $ f(x^{n+m})\in Z(\mathcal{B})$ for all $ x\in\mathcal{B}.$ And if we replace $x$ by $f^{-1}(x)$, we obtain 
$ x^{n+m}\in Z(\mathcal{B})$ for all $ x\in\mathcal{B}$ and  $\mathcal{B}$ is commutative (according to \cite{Moumen}).
\end{proof}
\begin{app}\label{a1}
Let $\mathbb{R}$ be the  field of real numbers. $\mathcal{B}=\mathcal{\mathcal{M}}_{2}(\mathbb{R})$
endowed with usual's matrix addition and  multiplication and the norm $\Vert .\Vert_{1}$ defined by $ \Vert A \Vert_{1}=\displaystyle\sum_{\substack{1\leq i,j \leq 2}}\vert a_{i,j}\vert $ for all $ A=(a_{i,j})_{1\leqslant i,j \leqslant2}\in\mathcal{B}$.
Let $f$  the automorphism of $\mathcal{B}$ defined by $f(M)=AMA^{-1}$ where 
 $A=\begin{pmatrix} 2 & 1 \\
1 & 2\end{pmatrix}$ ($f$ is a continuous automorphism).
Let $\mathcal{H}$ be a nonempty open subset of $\mathcal{B}$ included in $Z(\mathcal{B}).$
For all $ (A,B) \in \mathcal{H} \times \mathcal{H}  $ and for all $(p,q)\in\mathbb{N}^{*} \times \mathbb{N}^{*} $ we have
$ A^{p}\in Z(\mathcal{B})$ and $ B^{q}\in Z(\mathcal{B}),$ then  $f(A^{p}B^{q})+ A^{p}\circ B^{q}\in Z(\mathcal{B}).$
By Theorem $\ref{a},$ we conclude that $\mathcal{B}$ is commutative.
Consequently $\mathcal{H}=\emptyset$.\\
We conclude that  the only open subset included in $Z(\mathcal{B})$ is the empty set.
\end{app}
\begin{app}\label{a2}
Let $E$ be a normed space over $\mathbb{K}$ ($\mathbb{R}$ or $\mathbb{C}$). The space $\mathcal{B}={\displaystyle{\mathcal L_{c}}}(E)$  of continuous linear applications from $E$ to $E$,
endowed with  addition and  composition and the norm  defined by
$\parallel T \parallel =\displaystyle\sup_{\parallel x \parallel \leq 1}  \parallel T(x)\parallel$ for all
$T\in \mathcal{B} $, is a normed algebra over $\mathbb{K}$.
Let ${G}$ be the subspace of $\mathcal{B}$ defined  by ${G}=\{ \lambda I_{E} \,|\,\lambda\in\mathbb K \} $, where
$I_{E}$ is the identity of $E$, we observe that $G\subset Z(\mathcal{B}),$ according to Application $\ref{a1}$ we conclude that  the interior of $G$ is empty  because $\mathcal B $ is not commutative and $Int(G)\subset Int(Z(\mathcal{B}))=\emptyset$.
\end{app}
The following theorem can be shown using  the same techniques as the proof of  Theorem $1$.
%It is simple to establish the following theorem using a proof identical to that one with a few minor variations.
\begin{theorem}\label{d}
Let $\mathcal{B}$ be a prime Banach algebra, $\mathcal{H}_{1}$ and $\mathcal{H}_{2}$   two non-void open subsets of $\mathcal{B}.$
If $\mathcal{B}$ admits a  continuous automorphism $ f $  satisfying one of the following conditions\\
$i.\:\:(\forall (x,y)\in \mathcal{H}_{1}\times \mathcal{H}_{2})(\exists (p,q)\in\mathbb N^{*}\times \mathbb N^{*})  \;such\; that\;f(x^{p}y^{q})-x^{p} \circ y^{q}\in Z(\mathcal{B}).$\\
 $ii.\:\:(\forall (x,y)\in \mathcal{H}_{1}\times \mathcal{H}_{2})(\exists (p,q)\in\mathbb N^{*}\times \mathbb N^{*})  \;such\; that\;f(x^{p} \circ y^{q})-[x^{p} , y^{q}]\in Z(\mathcal{B}).$\\
Then $ \mathcal{B} $ is commutative.
\end{theorem}
\begin{app}
Let $n\geq 3$ be an integer and $\mathcal{B}=\mathcal{\mathcal{M}}_{n}(\mathbb{K})$ where $\mathbb{K}(=\mathbb{R}\; \mbox{or}\;= \mathbb{C}$). Let  $\mathcal{A}$ the set of all $n\times n$ strictly upper  triangular real or complexes matrices. Endowed with the usual operations on matrices and the norm $\|.\|_1$ defined by 
$\|A\|_1 = \displaystyle\sum_{1\leq i,j\leq n}|a_{ij}|$ for all $A=(a_{i,j})_{1\leq i,j\leq n}\in\mathcal{B}.$\\  $\mathcal{B}$ is a non-commutative real prime Banach algebra.
Observe that 
\begin{center}
     $f(x^{n}y^{n})-[x^{n},y^{n}]=0\in Z(\mathcal{B}) \; \; \forall (x,y)\in \mathcal B^2 $
\end{center}
where $f$ is the continuous automorphism of $\mathcal{B}$ defined by $f(M)=P^{-1}MP$ and $P\in\mathcal{B}$ is an  invertible matrix. Hence, it follows from Theorem \ref{a} that $\mathcal{A}$ is not a open subset in $\mathcal{B}$ .
\end{app}
%We close this article with the following corollary.
\begin{corollary}\label{d1}
Let $\mathcal{B}$ be a prime Banach algebra and $D$  a dense part of  $\mathcal{B}.$
 If $\mathcal{B}$ admits a  continuous automorphism $ f $ such that
$$
\exists (p,q)\in\mathbb N^{*}\times \mathbb N^{*}: f(x^{p} \circ  y^{q})+[x^{p},y^{q}]\in Z(\mathcal{B})\,\,\,\mbox{for all}\,\,\,x,y\in D,
$$
then $\mathcal{B}$ is commutative.
\end{corollary}
\begin{proof}
Let $x,y\in\mathcal{B},$ there exist two  sequences $(x_{k})_{k\in\mathbb{N}}$ and $(y_{k})_{k\in\mathbb{N}}$  in $ D $ converging to $x$ and $y$.
Since $(x_{k})_{k\in\mathbb{N}}\subset D $ and $(y_{k})_{k\in\mathbb{N}}\subset D $, so
\begin{center}
$f((x_{k})^{p} \circ (y_{k})^{q})+[(x_{k})^{p},(y_{k})^{q}]\in Z(\mathcal{A})$ for all $k\in\mathbb{N}.$
\end{center}
Since  $ d $ is continuous then  the sequence  $(f((x_{k})^{p} \circ (y_{k})^{q})+[(x_{k})^{p},(y_{k})^{q}])_{k\in\mathbb{N}}$ converges to
$ f(x^{p} \circ  y^{q})+[x^{p},y^{q}] $, knowing that $Z(\mathcal{B})$ is closed, then $ f(x^{p}\circ y^{q})+[x^{p},y^{q}]\in Z(\mathcal{B}).$\\
We conclude that :
$$ \exists (p,q)\in\mathbb N^{*}\times \mathbb N^{*}\;: f(x^{p} \circ  y^{q})+[x^{p},y^{q}]\in Z(\mathcal{B})\,\,\,\mbox{for all}\,\,\,x,y\in\mathcal{B}.$$ By Theorem $\ref{b},$ we get the desired conclusion.
\end{proof}
%\begin{remark}
%If we replace a dense part of $\mathcal{B}$ for any of the open sets in the preceding theorems, we can get conclusions that are similar to those in Corollary $\ref{d1}$.
%\end{remark}
\section{Generalized derivations in Banach algebras}
In this section, we examine some analogous conclusions for generalized derivations in Banach algebras.

\begin{theorem}\label{ga}
Let $\mathcal{B}$ be a  prime Banach algebra, $\mathcal{H}_{1}$ and  $\mathcal{H}_{2}$ non-void open subsets of $\mathcal{B}$. If $\mathcal{B}$ admits  a continuous  generalised derivation $F$ associated with a non-injective derivation $d$  satisfying  
 for all $(x,y)\in \mathcal{H}_{1}\times \mathcal{H}_{2}$ there are strictly positive integers $p',q'$ such that 
  $\;F(x^{p'}y^{q'})+x^{p'}\circ y^{q'}\in Z(\mathcal{B}),$
then  $ \mathcal{B}$ is   commutative or 
$d(Z( \mathcal{B}))=\{0\}$.
\end{theorem}

\begin{proof}
If $F=0,$ then our hypothesis reduces to $x^{p'} \circ y^{q'}\in Z(\mathcal{B})\;$  for all  $(x,y)\in \mathcal{H}_{1}\times \mathcal{H}_{2}$ so the required result follows from Theorem $1$ \cite{Moumen}. Hence we may suppose that $F\neq 0.$\\
 For all $(p',q')\in\mathbb N^{*}\times \mathbb N^{*}$
 let us define the following sets :
$$ O_{p',q'}=\{(x,y)\in \mathcal{B}^{2} \mid F(x^{p'}y^{q'})+ x^{p'}\circ y^{q'}\notin Z(\mathcal{B} )\}\;\; \mbox{and} \;\;     F_{p',q'}=\{(x,y)\in \mathcal{B}^{2} \mid F(x^{p'}y^{q'})+ x^{p'}\circ y^{q'}\in Z(\mathcal{B} )\}.$$
Using a similar approach with Theorem \ref{a}, we arrive at (where $(n,m)\in\mathbb N^{*}\times \mathbb N^{*}$)
\begin{equation*}
  F(x^{n}y^{m})+x^{n}\circ y^{m}\in Z(\mathcal{B})\;\;\mbox{for all}\;\;(x,y)\in \mathcal B^{2}.  
\end{equation*}

 Changing $x$ by $x^{m}$ and $y$ by $y^{n}$, We conclude that
$$
F(x^{nm}y^{mn})+ x^{nm}\circ y^{mn}\in Z(\mathcal{B})\,\,\,\mbox{for all}\,\,\, x,y\in\mathcal{B}.
$$
By permuting $x$ and $y$, we find that
$$
F(y^{nm}x^{mn})+ y^{nm}\circ x^{mn}\in Z(\mathcal{B})\,\,\,\mbox{for all}\,\,\, x,y\in\mathcal{B}.
$$
Then $$(F(x^{nm}y^{mn})+ x^{nm}\circ y^{mn})-(F(y^{nm}x^{mn})+ y^{nm}\circ x^{mn})\in Z(\mathcal{B})\,\,\,\mbox{for all}\,\,\, x,y\in\mathcal{B}.$$
Since  $$y^{nm}\circ x^{mn}=x^{nm}\circ y^{mn}\,\,\,\mbox{for all}\,\,\, x,y\in\mathcal{B},$$
we obtain
$$
F([x^{nm},y^{nm}])\in Z(\mathcal{B})\,\,\,\mbox{for all}\,\,\, x,y\in\mathcal{B}.
$$
Let $x\in\mathcal B$ and  $a=x^{nm}.$  We have $P(t)=F([a,(y+ta)^{nm}])\in Z(\mathcal{B})$ 
for all $(y,t)\in\mathcal{B}\times \mathbb{R}$.
This can be written as
 $$ P(t)=F([a,(y+ta)^{mn}])=\displaystyle\sum_{k=0}^{mn}t^{k}F([a,A_{nm-k,k}]),$$
where $ A_{nm-k,k} $ denotes the sum of all terms in which $ y $ appears exactly $ nm-k $ times and $ a $
appears exactly $ k $ times. By  Lemma $\ref{4},$
we conclude that $F([a,A_{nm-k,k}(x,y)])\in Z(\mathcal{B})$ for all $ 0\leq k\leq nm.$ The coefficient of $ t $ in this polynomial is $F([a,A_{nm-1,1}])$ and $ A_{nm-1,1}=\displaystyle\sum_{k=0}^{nm-1}a^{nm-1-k}y a^{k}.$ Since
\begin{align*}
 \sum_{k=0}^{nm-1}[a,a^{nm-1-k}y a^{k}]&=[a,a^{nm-1}y ]+[a,a^{nm-2}y a^{}]+[a,a^{nm-3}y a^{2}]+...+[a,y a^{nm-1}]\\
 &=(aa^{nm-1}y-a^{nm-1}ya)+(aa^{nm-2}ya-a^{nm-2}yaa)+..+(ay a^{nm-1}-y a^{nm-1}a)\\
 &=(a^{nm}y-a^{nm-1}y.a)+(a^{nm-1}ya-a^{nm-2}ya^{2})+...+(a.y a^{nm-1}-y a^{nm})\\
 &=a^{nm}y-y a^{nm}\\
 &=[a^{nm},y].
\end{align*}
Then
$$F([a, A_{nm-1,1}(x,y)])=\sum_{k=0}^{nm-1}F([a,x^{nm-1-k}y x^{k}])=F([a^{nm},y])\in Z(\mathcal{B}).$$
Thus
$$
F([x^{n^{2}m^{2}},y])\in Z(\mathcal{B}) \,\,\,\mbox{for all}\,\,\, (x,y)\in\mathcal{B}^{2}.
$$
Now, we prove that the restriction of $d $ on $ Z(\mathcal{B}) $ is not injective.\\
For this, we consider $a\in\mathcal{B}$ such that $d(a)=0$ and $b$ be a non-zero  element of $Z(\mathcal{A}).$  For all $t\in\mathbb{R},$ we have
$$
S(t)=F[(a+tb)^{p},y^{p}])\in Z(\mathcal{B})\,\,\,\mbox{for all}\,\,\,y\in\mathcal{B},
$$
where $p=n^{2}m^{2}$. By reason of $b^{k}\in Z(\mathcal{B})$ $(\forall k\in\mathbb{N^{*}}),$ we can write
$$
S(t)=\sum_{k=0}^{p}\binom {k}{p} t^{k}F([a^{p-k}b^{k},y^{p}])=\sum_{k=0}^{p}t^{k}A_{k}(a,b,y)\in Z(\mathcal{B})\,\,\,\mbox{for all}\,\,\,y\in\mathcal{B}
$$
where $A_{k}(a,b,y)=\binom {k}{p} F([a^{p-k}b^{k},y^{p}])$, according to Lemma $\ref{4},$ we conclude that $A_{k}(a,b,y)\in Z(\mathcal{B})$ for all $0\leq k\leq p$. In particular for $k=p-1$, we get $$p F([ab^{p-1},y^{p}])\in Z(\mathcal{B}),$$ we can  reduce by $p$ (because $Z(\mathcal{B})$ is a subspace of $ \mathcal{B}  )$ and we obtain $$F([ab^{p-1},y^{p}])\in Z(\mathcal{B}),$$ therefore $F(b^{p-1}[a,y^{p}])\in Z(\mathcal{B})$ (because $b^{p-1}$ is belongs to $Z(\mathcal{B})$). Since $$F(b^{p-1}[a,y^{p}])=b^{p-1}F([a,y^{p}])+[a,y^{p}]d(b^{p-1})\,\,\,\mbox{for all}\,\,\,y\in\mathcal{B}$$ and $b^{p-1}F([a,y^{p}])$ is an element of $Z(\mathcal{B}),$ then $d(b^{p-1})[a,y^{p}]\in Z(\mathcal{B})\,\,\,\mbox{for all}\,\,\,y\in\mathcal{B}.$ \\
As $d(b^{p-1})$ is a non-zero element of $Z(\mathcal{B}),$ by Fact \ref{rm}, we conclude that $$[a,y^{p}]\in Z(\mathcal{B})\,\,\,\mbox{for all}\,\,\,y\in\mathcal{B}.$$ Therefore, for all $t\in\mathbb{R}$ we have
$$Q(t)=[a,(a+ty)^{p}]\in Z(\mathcal B)\,\,\,\mbox{for all}\,\,\,y\in\mathcal{B}.$$
That is
$$
Q(t)=\sum_{k=0}^{p}t^{k}[a,B_{p-k,k}(a,y)] \,\,\,\mbox{for all}\,\,\,y\in\mathcal{B}
$$
where $ B_{p-k,k}(a,y) $ denotes the sum of all terms in which $ a $ appears exactly $ p-k $ times and $ y $
appears exactly $ k $ times. By Lemma $\ref{4},$
we deduce that $[a,B_{p-k,k}(a,y)]\in Z(\mathcal{B})$ for all $k\leq p.$ \\
The coefficient of $ t $ in this polynomial is $[a,B_{p-1,1}(a,y)]$ where $ B_{p-1,1}(a,y)=\sum_{k=0}^{p-1}a^{p-1-k}y a^{k},$ therefore
$[a, B_{p-1,1}(a,y)]=\sum_{k=0}^{p-1}[a,a^{p-1-k}y a^{k}]\in Z(\mathcal{B}),$
we observe that
\begin{align*}
 \sum_{k=0}^{p-1}[a,a^{p-1-k}y a^{k}]&=[a,a^{p-1}y ]+[a,a^{p-2}y a^{}]+[a,a^{p-3}y a^{2}]+...+[a,y a^{p-1}]\\
 &=aa^{p-1}y-a^{p-1}ya+aa^{p-2}ya-a^{p-2}yaa+...+ay a^{p-1}-y a^{p-1}a\\
 &=a^{p}y-a^{p-1}ya+a^{p-1}ya-a^{p-2}ya^{2}+...+ay a^{p-1}-y a^{p}\\
 &=a^{p}y-y a^{p}\\
 &=[a^{p},y].
\end{align*}
Consequently,
$$
[a^{p},y]\in Z(\mathcal{B})\:\:\mbox {for all}\:\: y\in\mathcal{B}.
$$
There are $2$ cases:\\
\begin{description}
\item[Case 1] : If $a^{p}\in Z(\mathcal{B}),$ since $d$ is injective in $Z(\mathcal{B})$ and $d(a^{p})=0$ then  $a^{p}=0,$  therefore $a=0$ and $d$ is injective.\\
\item[Case 2]: If $a^{p}\notin Z(\mathcal{B})$ then $d(\mathcal{B})\subset \mathcal{B}$ where $d$ is the inner derivation associated to $a^{p},$  by Fact $\ref{rm}$ we have $\mathcal{B}\subset Z(\mathcal{B}),$ so $d$ is injective.\\
\end{description}
Lastly, we conclude that the restriction of $ d $ on $ Z(\mathcal{B}) $ is non-injective, then there is a nonzero element $a$  of $Z(\mathcal{B})$  such that $ d(a)=0,$
we have
$ F([(x+ta)^{p},y] \in Z(\mathcal{B})$ for all $(x,y)\in\mathcal{B}$ and $ t\in \mathbb{R} .$
Since $a^{k}\in Z(\mathcal{B})\,\,\,\mbox{for all}\,\,\,$ $k\in\mathbb{N^{*}},$ we can write
             $$ F([(x+ta)^{p},y])=\sum_{k=0}^{p} \binom {p}{k}t^{k}F([ a^{p-k}x^{k},y])\in Z(\mathcal{B}),$$
 according to  Lemma $\ref{4}$, we conclude  that $ \binom {p}{k}F([a^{p-k}x^{k},y])\in Z(\mathcal{B})$ for all $0\leq k \leq p.$
In particular for $k=1,$ we have $ pF([a^{p-1} x,y])\in Z(\mathcal{B})$ and $F([a^{p-1} x,y])\in Z(\mathcal{B}),$ since $a^{p-1}\in Z(\mathcal{B})$  we have $F(a^{p-1}[x,y])\in Z(\mathcal{B})$, then $a^{p-1}F([x,y])\in Z(\mathcal{B})$ (because $ d(a^{p-1})=0 $). So $a^{p-1}\in Z(\mathcal{B})\smallsetminus\{0\},$
according to  Fact \ref{rm} we have
\begin{center}
$F([x,y])\in Z(\mathcal{B})  \; \; \mbox{for all} \; \; (x,y)\in\mathcal{B}.$
\end{center}

If $d(Z( \mathcal{B}))\neq \{0\}$, then, there is $z\in Z( \mathcal{B})$ such that $d(z)\neq 0$ and we have 
$$F([x,yz])\in Z( \mathcal{B}) \; \; \mbox{for all} \; \; (x,y)\in\mathcal{B}.$$
That is 
$$F([x,y]z)\in Z( \mathcal{B}) \; \; \mbox{for all} \; \; (x,y)\in\mathcal{B}.$$
And 
$$F([x,y])z+[x,y]d(z)\in Z( \mathcal{B}) \; \; \mbox{for all} \; \; (x,y)\in\mathcal{B}.$$
Since 
$$F([x,y])z\in Z( \mathcal{B}) $$
therfore 
$$[x,y]d(z)\in Z( \mathcal{B}) \; \; \mbox{for all} \; \; (x,y)\in\mathcal{B}$$
by Fact \ref{rm} we conclude that $$[x,y]\in Z( \mathcal{B}) \; \; \mbox{for all} \; \; (x,y)\in\mathcal{B}.$$
Replacing $y$ be $xy$, we abtain 
 $$x[x,y]\in Z( \mathcal{B}) \; \; \mbox{for all} \; \; (x,y)\in\mathcal{B}.$$
Then $x\in Z( \mathcal{B})$ or $[x,y]=0$,  both shows that $\mathcal {B}$ is commutative.
\end{proof}
We show the following result by traveling on the same path.

\begin{theorem}\label{145}
Consider $\mathcal H_1$ and $\mathcal H_2$ as nonvoid open subsets of  a prime Banach algebra $\mathcal{B}$. If $\mathcal{B}$ admits  a continuous  generalised derivation $F$ associated with a non-injective derivation $d$  satisfying  
 for all $(x,y)\in \mathcal{H}_{1}\times \mathcal{H}_{2}$ there are strictly positive integers $p$ and $q$ such that 
  $\;F(x^{p}\circ y^{q})+[x^{p}, y^{q}]\in Z(\mathcal{B}),$
then  $ \mathcal{B}$ is   commutative or 
$d(Z( \mathcal{B}))=\{0\}$.
\end{theorem}

\section{Some examples}
The example that follows demonstrates that our theorems' assertion that " $\mathcal{H}_{1}$ and $\mathcal{H}_{2}$ are open" is not superfluous.
\begin{example}
%Let $\mathbb{R}$ be the  field of real numbers. 
Let $\mathcal{B}=\mathcal{\mathcal{M}}_{2}(\mathbb{R})$ the ring of square matrices with real inputs of order $2 \times 2$
endowed with usual's matrix addition and  multiplication and the norm  defined by \\$ \Vert A \Vert_{1}=\displaystyle\sum_{\substack{1\leq i,j \leq 2}}\vert a_{i,j}\vert $ for all $ A=(a_{i,j})_{1\leqslant i,j \leqslant2}\in\mathcal{B},$  is a real  prime  Banach algebra.\\
Let
$ \mathcal{F}_{1}=\Big\{\begin{pmatrix} t & 0 \\
0 & t\end{pmatrix} \,\,|\,\, t \in \mathbb{R}  \Big\}$ and
$ \mathcal{F}_{2}=\Big\{\begin{pmatrix} t & 0 \\
0 & t\end{pmatrix} \,\,|\,\,t \in \mathbb{R^{+}}  \Big\}.$\\
$\mathcal{F}_{1}$ is not open in $\mathcal{B},$ indeed, we have to show that the complement of $\mathcal{F}_{1}$ is not closed. For this, we consider the sequence
$\Big(\begin{pmatrix} 1+\dfrac{1}{n} & \dfrac{-1}{n} \\
\dfrac{1}{n} & 1+\dfrac{1}{n}\end{pmatrix} \Big)_{n \in \mathbb{N^{*}}}$ in $\mathcal{F}_{1}^{c}$ complement of $\mathcal F_{1}$
 who converge to $\begin{pmatrix} 1 & 0 \\
0 & 1\end{pmatrix}\notin\mathcal{F}_{1}^{c} ,$ %therefore $(t_{n})_{n\in\mathbb{N}}$ converge to $t\in\mathbb{R},$ then
then $\mathcal{F}_{1}^{c}$ is not closed, that is $\mathcal{F}_{1}$ is not open in $\mathcal{B}.$\\
%Let $id$ the identity mapping on   $\mathcal{B}$. 
Let $ A=\begin{pmatrix} a & 0 \\
0 & a\end{pmatrix} \in\mathcal{F}_{1}$ and $ B=\begin{pmatrix} b & 0 \\
0 & b\end{pmatrix} \in\mathcal{F}_{2}$. For all  $ (p,q)\in\mathbb{N}^{2},$ we  have
  $$A^{p}=\begin{pmatrix} a^{p} & 0 \\
0 & a^{p}\end{pmatrix} \;   and \; B^{q}=\begin{pmatrix} b^{q} & 0 \\
0 & b^{q}\end{pmatrix}.$$
So $$ A^{p}B^{q}=\begin{pmatrix} a^{p}b^{q} & 0 \\
0 & a^{p}b^{q}\end{pmatrix}
,\;
A^{p}\circ B^{q}=\begin{pmatrix} 2a^{p}b^{q} & 0 \\
0 & 2a^{p}b^{q}\end{pmatrix}\;
and \; [A^{p},B^{q}]=\begin{pmatrix} 0 & 0 \\
0 & 0\end{pmatrix}.$$
We conclude that
$$Id(A^{p}\circ B^{q})=Id([A^{p},B^{q}])=Id(A^{p}B^{q})=\begin{pmatrix} 0 & 0 \\
0 & 0\end{pmatrix}.$$
(where $Id$ is the identity application of $\mathcal{B}$). We obtain:\\
1. $Id(A^{p}B^{q})+A^{p}\circ B^{q}\in Z(\mathcal{B})$ and $Id(A^{p}B^{q})-A^{p}\circ B^{q}\in Z(\mathcal{B}),$ \\
2. $Id(A^{p}\circ B^{q})+[A^{p},B^{q}]\in Z(\mathcal{B})$ and $Id(A^{p}\circ B^{q})-[A^{p},B^{q}]\in Z(\mathcal{B}).$  \\
But $\mathcal{B}$ is not commutative.
\end{example}
The example that follows demonstrates that we cannot change $\mathbb{R}$ or $\mathbb{C}$ by $\mathbb{F}_{3}=\mathbb{Z}/3\mathbb{Z}$ in the statement of previous theorems.
\begin{example}
Let $\mathcal{B}=(\mathcal{M}_{2}(\mathbb{Z}/3\mathbb{Z}),+,\times,.)$ the 2-torsion free prime  Banach algebra of square matrices of size $ 2 $ at coefficients in $ \mathbb{Z}/3\mathbb{Z} $ with usual matrix addition and matrix multiplication.\\ The norm  is defined by $ \Vert A \Vert_{1}=\displaystyle\sum_{\substack{1\leq i,j \leq 2}}\vert a_{i,j}\vert $ for all $ A=(a_{i,j})_{1\leqslant i,j \leqslant2}\in\mathcal{B}$ with $\vert . \vert$  is the norm defined on $\mathbb{Z}/3\mathbb{Z}$ by:
\begin{center}
$\vert \overline{0}\vert =0 $, $\vert \overline{1}\vert =1 $ and $\vert \overline{2}\vert =2 .$
\end{center}
Observe that
$\mathcal H=\{\begin{pmatrix} a & 0 \\
0 & a \end{pmatrix} \:|\: a \in \mathbb{Z}/3\mathbb{Z} \}$ is open in  $\mathcal{B}$. Indeed, let $A \in \mathcal H$  the open ball $B (A,1)=\{X\in \mathcal{B}$ such that $ \Vert A-X \Vert_{\infty}< 1 \}=\{A\} \subset \mathcal H,$ %so $ H $ is  neighborhood  of $ A $
then $\mathcal H $ is a non-void open subset of $\mathcal{B}.$
For  all  $(p,q)\in\mathbb N^* \times \mathbb N^*$ and for all $(A,B) \in\mathcal H \times \mathcal H$ we have :\\
%\begin{center}
 1. $Id(A^{p}B^{q})+A^{p}\circ B^{q}\in Z(\mathcal{B})$ and $Id(A^{p}B^{q})-A^{p}\circ B^{q}\in Z(\mathcal{B}).$  \\
 2. $Id(A^{p}\circ B^{q})+[A^{p},B^{q}]\in Z(\mathcal{B})$ and $Id(A^{p}\circ B^{q})-[A^{p},B^{q}]\in Z(\mathcal{B}).$  \\
 So all the conditions are verified except $\mathbb{K}=\mathbb{R}$ or $\mathbb{K}=\mathbb{C}$, but $ \mathcal{B}$ is not commutative.
\end{example}

\end{document}